\documentclass[12pt]{iopart}
\usepackage{iopams}

\expandafter\let\csname equation*\endcsname\relax  
\expandafter\let\csname endequation*\endcsname\relax
\usepackage{amsmath}
\usepackage{graphicx}
\usepackage{bm}
\usepackage[usenames, dvipsnames]{color}
\usepackage{booktabs}
\usepackage{multirow}
\usepackage[caption=false]{subfig}
\usepackage[font=footnotesize,labelfont=bf, justification=justified, format=plain]{caption}
\usepackage{hyperref}  
\usepackage{cleveref}  

\newcommand\Rey{\mathrm{Re}}

\begin{document}
    \title{Tracking complex singularities of fluids on log-lattices}
    \author{Quentin Pikeroen$^1$, Amaury Barral$^1$, Guillaume Costa$^1$, Ciro Campolina$^{2}$, Alexei Mailybaev$^3$ \& Berengere Dubrulle$^1$}
    \address{
        $^1$SPEC/IRAMIS/DSM, CEA, CNRS, University Paris-Saclay, CEA Saclay, 91191 Gif-sur-Yvette, France
    }
	\address{$^2$Universit\'{e} C\^{o}te d’Azur, CNRS, Inria, Institut de Physique de Nice, France}
	\address{$^3$Instituto de Matem\'{a}tica Pura e Aplicada, 22460-320 Rio de Janeiro, Brazil}
    \ead{berengere.dubrulle@cea.fr}

    \begin{abstract}
 In 1981, Frisch and Morf~\cite{FM81} postulated the existence of complex singularities in solutions of Navier-Stokes equations.
        Present progress on this conjecture is hindered by the computational burden involved in simulations of the Euler equations or the Navier-Stokes equations at high Reynolds numbers.
        We investigate this conjecture in the case of fluid dynamics on log-lattices, where the computational burden is logarithmic concerning ordinary fluid simulations.
        We analyze properties of potential complex singularities in both 1D and 3D models for lattices of different spacings.
        Dominant complex singularities are tracked using the singularity strip method to obtain new scalings regarding the approach to the real axis and the influence of normal, hypo and hyper dissipation.
    \end{abstract}

    \maketitle

    \section{Introduction}
 Viscous fluids dissipate mechanical energy into heat due to the first law of thermodynamics.
    Observations and numerical simulations reveal that this dissipation is not homogeneous within the flow but occurs via spatially or temporally intermittent bursts, a phenomenon classically referred to as intermittency.
    Moreover, after spatial and temporal averaging, the mean energy dissipation becomes independent of the viscosity in the inviscid limit, according to the empirical ``zeroth law of turbulence''.
    Onsager explained these observations in 1949~\cite{[O49]}, conjecturing that strong enough singularities in the inviscid flow could provide an anomalous dissipation.
    While this conjecture has been proven mathematically~\cite{Isett18}, its application to fluids is still debated.
	In fact, the developement of finite-time singularities in Euler flows is until now an unsolved problem~\cite{gibbon2008three},
	while the same question formulated for the Navier-Stokes equations is among the open Millennium Prize Problems of the Clay Mathematics Institute~\cite{fefferman2006existence}.
    This debate concerns the existence of singularities in real space.
    In 1981, Frisch and Morf~\cite{FM81} paved the way to another possibility based on the existence of complex singularities.
    They proved on a simple one-dimensional non-linear Langevin system that the dynamics of such complex singularities could be directly connected to intermittency, as dissipation bursts occur whenever a complex singularity approaches the real axis.

 Since then, this scenario was also confirmed in the one-dimensional Burgers equation -- a 1D surrogate of the Navier-Stokes equation.
    In this system, real singularities can be observed in the inviscid limit and manifest as shocks, i.e.\ finite jumps in the velocity.
    Shocks dissipate energy in agreement with the dissipation anomaly~\cite{EyinkLN}.
    They correspond to the collapse of two complex conjugate singularities onto the real axis~\cite{FF83, SRE96}.
    When a viscosity $\nu$ is added, the singularities are repelled from the real axis, the closest one being constantly at a distance greater than $O(\nu^{3/4})$ to the real axis.
    The complex singularities follow Calogero-Moser (CM) dynamics~\cite{Ca78}, with long-range interactions (decaying in $1/r$).
    There is an exact mapping between such CM dynamics and the solution of the PDE, which can be described exactly via pole decomposition coupled to the integration of the CM equations~\cite{SRE96}.

 The generalization of these findings to 3D is challenging~\cite{frisch2003singularities}. The computational burden to resolve the Navier-Stokes  equation for a fluid with typical velocity $U$ and length $L$ scales like $\Rey^{3}$, where $\Rey\sim UL/\nu$ is the Reynolds number.
    Most of the earlier attempts to track complex singularities in the inviscid limit were performed using the ``singularity strip'' method~\cite{sulem1983tracing}, which is based on the observation that the behaviour of the energy spectrum at large wavenumber $k$ is dominated by the position of the singularity closest to the real axis, and decays like
    $\exp(-2\delta k)$, where $\delta$ is the imaginary part of corresponding singularity.
    Fitting the large wavenumber tail of the energy spectrum as a function of time, one then gets an estimate of
    $\delta(t)$, and a real singularity occurs when $\delta(t)=0$.
    So far, studies have only identified exponentially decaying regimes for $\delta(t)$~\cite{siegel2009calculation} which suggests the absence of finite time blow-up.
    However, we cannot guarantee that this extrapolation is correct due to numerical limitations.

 New perspectives on these issues were opened recently by Campolina and Mailybaev~\cite{CM21}, exploring fluid dynamics on log-lattices.
    This technique may be viewed as a generalization of the so-called ``shell models''~\cite{gloaguen85,biferale03} and solves the equations of motion in Fourier space using a sparse set of Fourier modes.
    The modes are evenly spaced points in log space (``logarithmic lattices''). They interact via nonlinear equations derived from the fluid equations by substituting for the convolution product a new operator, which can be seen as a convolution on the log-lattice, while preserving most symmetries of the original equation.
    The model is valid for all dimensions.
    In 1D, it was shown to encompass~\cite{CM21} the dyaic and Sabra shell models of turbulence~\cite{gloaguen85,biferale03}.
    In 3D, its solutions have the same behaviour as the Navier-Stokes equation in Fourier space (energy spectrum, energy transfers), over an unprecedented wide range of scales~\cite{CM21}.
    In the inviscid equations, a finite-time blow-up is observed~\cite{CM18} in connection with a chaotic attractor that propagates at a constant average speed in a renormalized Fourier space, like a wave.
    However, Campolina and Mailybaev did not attempt to track possible complex singularities in connection with such a blow-up.

   This is the purpose of the present paper.
    In the first part, we validate the close connection between fluid dynamics on log-lattice and real fluid dynamics by focusing on the 1D Burgers equation, where dominant complex singularities are tracked using the singularity strip method.
    In the second part, we extend this technique to 3D to obtain new scalings regarding the approach to the real axis and the influence of normal, hypo and hyper dissipation.

    \section{Log-lattice framework}

    \subsection{Definitions and notations}
    We consider a $d$-dimensional complex vector field $u(t,k) = (u_1,\dots,u_d)$ depending on time $t \in \mathbb{R}$ and on the wave vector $k = (k_1,\dots,k_d)$.
    We shall interpret $u$ as the Fourier components of the velocity field.
    For this reason, we require them to satisfy the Hermitian symmetry $u(t,-k) = \overline{u(t,k)}$ with respect to $k$, which is the Fourier property of a real-valued function in physical space.
    The wave vector $k$ is embedded on a \textit{logarithmic lattice} (in short, \textit{log-lattice}), which means that its components follow geometric progressions $k = k_0(\pm\lambda^{m_1},\dots,\pm\lambda^{m_d})$ for integers $m_1,\dots,m_d$, where $k_0 = 2\pi$ is a fixed positive reference scale, and $\lambda>1$ is the spacing factor of the lattice.
    The dependence of $u$ on $t$ and $k$ is henceforth implicit and specified only when ambiguity prevails.

    \textit{Fluid dynamics on log-lattice}~\cite{CM21} is the set of vector fields $u$ which are solutions of the equations
    \begin{subequations}
        \label{eq:loglatice}
        \begin{align}
            k_\beta u_\beta&=0,\label{eq:loglattice_incomp}\\
            \partial_t u_\alpha+ik_\beta (u_\alpha* u_\beta)&=-i k_\alpha p-\nu k^{2\gamma} u_\alpha+f_\alpha,\label{eq:loglattice_mom}\\
            (u_\alpha * u_\beta)(k)&=\sum_{q+r=k} u_\alpha(q) u_\beta(r)\label{eq:loglattice_conv},
        \end{align}
    \end{subequations}
    where $p$ is the complex pressure field that enforces incompressibility~\eqref{eq:loglattice_incomp}, $f$ is a vectorial forcing, and $\nu$ is a non-negative viscosity parameter.
    When $\nu > 0$, the exponent $\gamma$ measures the dissipation degree: we say the flow has \textit{viscous} (or \textit{usual}) \textit{dissipation} if $\gamma = 1$, it has \textit{hypo-dissipation} if $\gamma<1$, and it has \textit{hyper-dissipation} if $\gamma>1$.
    Similarly to the dynamics of continuous media, system~\eqref{eq:loglatice} is the \textit{incompressible Navier-Stokes equations} on the log-lattice.
    When $\nu = 0$, the flow is \textit{inviscid}, and the system reduces to the \textit{incompressible Euler equations} on the log-lattice.

    The convolution in \cref{eq:loglattice_conv} defines triadic interactions on the logarithmic lattice, which are nontrivial only if the equation $\lambda^m=\pm\lambda^q\pm\lambda^r$ has integer solutions $m,q,r$.
    As shown in~\cite{CM21}, this is possible only for particular values of $\lambda$, which determine the number of possible interactions on the grid.
    In this paper, we consider the following three values: $\lambda=2$, with $3$ interactions per direction; $\lambda= \phi\approx 1.618$ (the golden number), with $6$ interactions per direction; and $\lambda=\sigma\approx 1.325$ (the plastic number), with $12$ interactions per direction.
    As $\lambda$ decreases from $2$ to $\sigma$, the density of nodes and the number of interactions on the grid increase.
    We recall, however, that the interactions for these log-lattices are all local.

    \subsection{Global quantities}
    By analogy with the Fourier representation of classical fluid flows, we define the global quantities representing the \textit{total energy} $E$ and the \textit{helicity} $H$ as
    \begin{align}
        E&=\sum_k \vert u\vert^2,\label{eq:total_energy}\\
        H&=\sum_k u_\alpha \overline{\omega}_\alpha\label{eq:helicity},
    \end{align}
    where $\omega_\alpha = \epsilon_{\alpha \beta \gamma} ik_\beta u_\gamma$ is the \textit{vorticity} field; here, $\epsilon_{\alpha \beta \gamma}$ is the Levi-Civita symbol.
    Regular solutions of the unforced three-dimensional inviscid system~\eqref{eq:loglatice} conserve these quantities in time~\cite{CM21}.

    Moreover, we define the energy spectrum $E(k)$ as
    \begin{equation}
        E(k)= \langle\vert u\vert^2\rangle_{S_k},
        \label{eq:defispecm}
    \end{equation}
    where the average $\langle\ \cdot\ \rangle_{S_k}$ is taken over the wave vectors in the shell $S_k$ delimited by spheres of radii $k$ and $\lambda k$.
    More explicitly,
    \begin{equation}
        \langle\vert u\vert^2\rangle_{S_k}=\frac{1}{N_k(\lambda k - k)}\sum_{k\le\vert q \vert < \lambda k } \vert u(q)\vert^2,
        \label{eq:defispec}
    \end{equation}
    where $N_k\sim (\log k)^{d-1}$ is the number of wave vectors in the shell $S_k$.

    \subsection{Regularity}
   The solutions of fluid dynamics equations on log-lattices~\eqref{eq:loglatice} share some regularity properties with the original models.
    The main mathematical results are for the inviscid Euler equations~\cite{CM21}.
    For this system, the local-in-time existence of strong solutions and a Beale-Kato-Majda blow-up criterion were proved.
    Exploiting the conservation of enstrophy, one proves the global regularity of two-dimensional flows.
    In the three-dimensional case, high-resolution numerical log-lattice simulations disclosed a finite-time blow-up, characterized by a chaotic wave travelling with constant average speed along a renormalized set of variables~\cite{CM18}.
    Such blow-up scenario was confirmed for $\lambda = \phi$ and $\lambda = \sigma$, presenting the same asymptotic blow-up scalings~\cite{CM21}.
    In the viscous case, numerical simulations suggest the expected global regularity of solutions.

    \subsection{Singularity strip method for log-lattices}
    If a potential singularity is due to an imaginary pole crossing the real axis, one can track its distance to the real axis via the \textit{singularity strip method}~\cite{sulem1983tracing}.
    This method considers the analytic continuation $u(z)$ of the physical-space velocity field and is based on the following property: if
    \begin{equation}
        u(z)\sim 1/(z-z_*)^\xi, \quad \text{for} \ z \to z_*
        \label{eq:singustripPole}
    \end{equation}
    in a neighborhood of the complex singularity $z_*=a+i\delta$, then its Fourier transform $\hat{u}_k$ satisfies
    \begin{equation}
        \hat{u}_k\sim k^{-d-\xi} e^{ika}e^{-\delta k}, \quad \text{as} \ k \to \infty.
        \label{eq:singustripLL}
    \end{equation}
    Asymptotics of~\eqref{eq:singustripLL} provide the corresponding exponential decay $E(k) \sim e^{-2\delta k}$ for the energy spectrum over a typical length $2\delta$.
    Therefore, one can measure the distance of the dominant pole to the real axis by monitoring the decay of the energy spectrum in Fourier space.
    A finite-time singularity at instant $t_b$ would occur if $\delta \to 0$ as $t \to t_b$.

    Extension of this notion to the log-lattice framework is natural. It relies on the observation that if a flow~\eqref{eq:loglatice} on log-lattice satisfies $u(k)\sim k^{-d-\xi} e^{-\delta k}$, then its inverse Fourier transform obeys a relation similar to~\eqref{eq:singustripPole}.
    Therefore, we can generalize the singularity strip method to log-lattices, where $2\delta$ is estimated from the slope of $\log E(k)$ as a function of $k$.


    \subsection{Numerical methods}

	Equations~\eqref{eq:loglatice} are numerically integrated using a technique analogous to viscous splitting.
	Considering a time step $dt$, we obtain $u(t+dt)$ from $u(t)$ employing the following strategy.
	Using $u(t)$ as initial condition, we first solve the inviscid equation
	\begin{equation}
		\label{eq:nsinviscid}
		\partial_t{u_\alpha} =P_{\alpha\beta} \left[- i k_\sigma (u_\beta\ast u_\sigma) + f_\beta\right],
	\end{equation}
	where $P_{\alpha\beta} = \delta_{\alpha\beta} - \frac{k_\alpha k_\beta}{k^2}$ accounts for the pressure term under the incompressibility hypothesis~\eqref{eq:loglattice_incomp}.
	For that, we use an explicit 4th order Runge-Kutta method.
	This yields $u(t+dt)_{\nu = 0}$.
	Then, the viscosity is taken into account through $u(t+dt) = u(t+dt)_{\nu = 0}e^{-\nu k^{2\gamma}dt}$.
	In this whole process, we adapt dynamically the time step $dt$.
	

    \section{1D Burgers equation}

    Before going to the full three-dimensional Navier-Stokes system on log-lattices, we take an intermediate step by studying the simpler one-dimensional Burgers equation.
    For this system, several exact mathematical results are available.
    This allows us to probe the singularity strip method on log-lattices, by comparing our numerical computations with the exact expected results.

    The one-dimensional Burgers equation on log-lattices is obtained from system~\eqref{eq:loglatice} as follows.
    We consider a compressible pressureless flow on a one-dimensional log-lattice.
    Mathematically, this translates into setting $p = 0$ and dropping \cref{eq:loglattice_incomp} from the system, which reduces to
    \begin{subequations}
        \label{eq:burgers}
        \begin{align}
            \partial_t u + u * \partial_x u&= -\nu k^{2\gamma} u+f,\label{eq:burgers_mom}\\
            (u * \partial_x u)(k)&=\sum_{q+r=k} ir \, u(q) u(r)\label{eq:burgers_conv}.
        \end{align}
    \end{subequations}
    
    It was shown~\cite{campolina2019fluid} that, up to a prefactor in the convolution~\eqref{eq:burgers_conv}, the Burgers equation on log-lattices is equivalent to well-known shell models of turbulence for specific choices of parameters.
    Particularly, when $\lambda = 2$, system~\eqref{eq:burgers} (but with a factor 2 added in the convolution and restricting to imaginary solutions) is the dyadic model~\cite{desnyansky1974evolution}, while for $\lambda = \phi$ (but with a factor $-\phi^2$ added in the convolution) it is the Sabra model~\cite{lvov1998improved} in a three-dimensional parameter regime (second invariant is not sign defined).
    Because of this relation with shell models of turbulence, the Burgers equation on the one-dimensional log-lattice inherits several results concerning the regularity of its solutions, which we briefly review now.

    For the dyadic model ($\lambda = 2$) with $\nu >0$, there are theorems~\cite{Ch08} for global existence of weak solutions (satisfying the energy inequality at almost all time), local regularity when $\gamma > 1/3$, global regularity when $\gamma \geq 1/2$, and finite-time blow-up when $\gamma < 1/3$ for sufficiently large initial conditions.
     Note that for the {\sl continuous} version of the 1D Burgers equations, global existence and analycity holds whenever $\gamma \geq 1/2$, while finite-time blow-ups are present whenever $\gamma<1/2$ \cite{Kiselev2008}.
    In contrast, we have  presently no rigorous statements about the dyadic model for the parameter range $1/3 \leq \gamma < 1/2$. This means that the mathematical techniques used in the currently available theorems are not sharp enough to separate the finite-time blow-up and the global regularity regimes.
    The finite-time singularity in the inviscid case was also rigorously established~\cite{katz2005finite}.

    For the viscous Sabra model ($\lambda = \phi$) with usual dissipation $\gamma = 1$, there are     
    proofs~\cite{constantin2006analytic} of global regularity of strong solutions.
    Like the Navier-Stokes equations, the dynamics of the Sabra model develops within finite degrees of freedom.
    Indeed, the finite dimensionality of the global attractor and the existence of a finite-dimensional inertial manifold were proved~\cite{constantin2006analytic}.
    On the other hand, the inviscid model has~\cite{CLT07} global-in-time existence of weak solutions with finite energy, local-in-time regularity, and a Beale-Kato-Majda blow-up criterion.
    Despite the absence of rigorous proofs, it is well-known~\cite{mailybaev2016spontaneously} that Sabra (in the three-dimensional parameter regime) develops a self-similar finite-time blow-up, characterized as a travelling wave in a renormalized system of variables (cf.~\cite{dombre1998intermittency}).
    Following the dynamical systems approach, such blow-up can be seen as a fixed-point attractor of the associated Poincar\'{e} map~\cite{mailybaev2013bifurcations}.

    To our knowledge, there are no systematic results about the development of singularities in Sabra with general dissipation exponents $\gamma$, nor in the case of our third lattice parameter $\lambda = \sigma$.

    \subsection{Inviscid flow}
\begin{figure}
        \subfloat[\label{fig:SpectBurgers}]{%
            \includegraphics[width=0.5\textwidth]{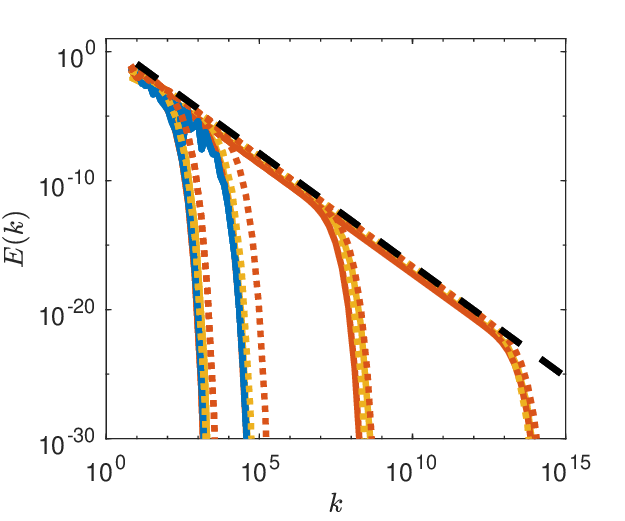}
        }
        \subfloat[\label{fig:OmaxAdimBurgers}]{%
            \includegraphics[width=0.5\textwidth]{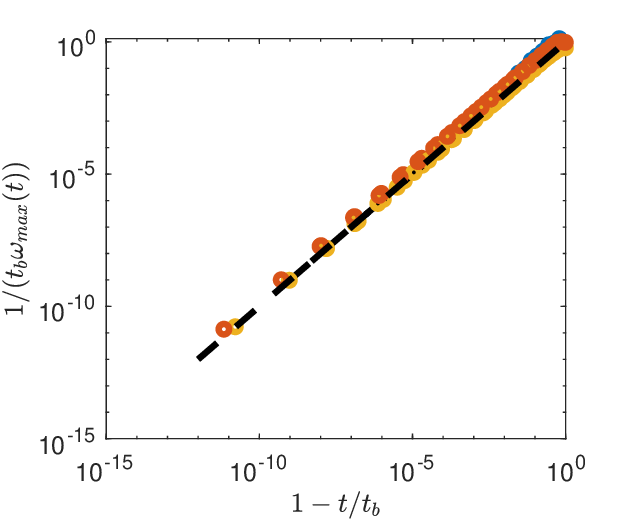}
        }\\
        \subfloat[\label{fig:deltaBurgers}]{%
            \includegraphics[width=0.5\textwidth]{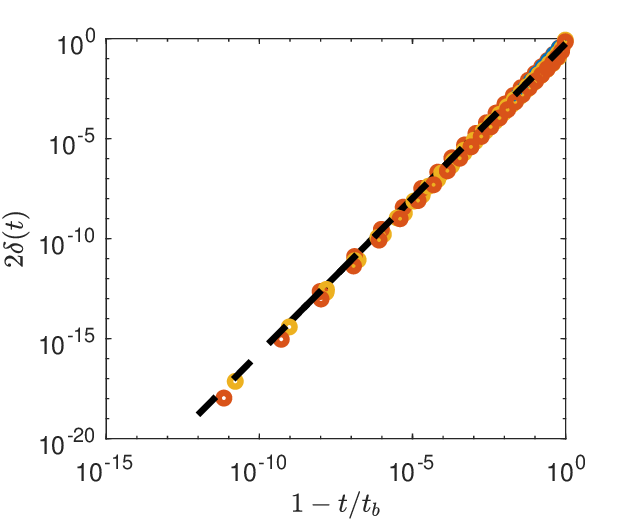}
        }
        \subfloat[\label{fig:deltaAdimBurgers}]{%
            \includegraphics[width=0.5\textwidth]{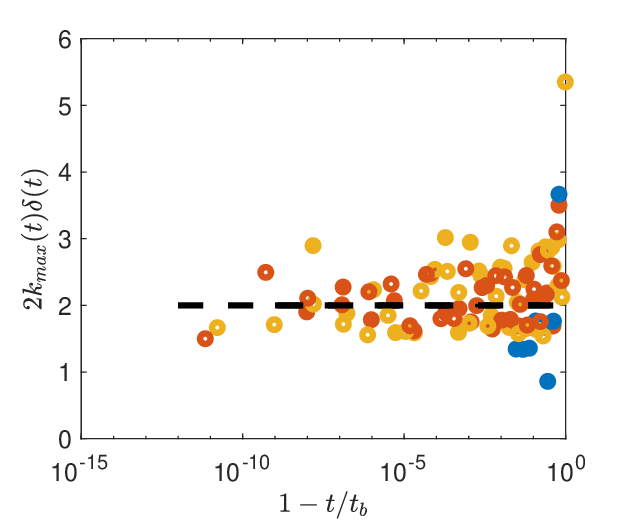}
        }
        ~\caption{\footnotesize Inviscid blow-up for the 1D Burgers equation for $\lambda=2$ (yellow), $\lambda=\phi$ (red) and $\lambda=\sigma$ (blue).
        Continuous lines and filled symbols indicate simulations with constant forcing, while dotted lines and open symbols indicate simulations without forcing.
            (\ref{fig:SpectBurgers})~Spectra at different renormalized relative time $\tau=1-t/t_b$. The black dotted line has a slope of $-1.733$.
            (\ref{fig:OmaxAdimBurgers})~Maximum value of the derivative $1/t_b \omega_{\max}$ as a function of $\tau$; The black dotted line is the theoretical value from \cref{eq:ssblowup}.
            (\ref{fig:deltaBurgers})~Width of the analyticity strip $2\delta$ as a function of $\tau$. The black dotted line has a slope given in \cref{tab:exp}.
            (\ref{fig:deltaAdimBurgers})~Renormalized width $k_{\max}\delta$ as a function of $\tau$. The black dotted line has a slope of $0$.}
        \label{fig:euler1d}
    \end{figure}

    We start with the inviscid ($\nu = 0$) Burgers equation with and without forcing.
    When forcing, initial conditions are equal to zero, and the forcing is equal to the imaginary unit $i$ on the first mode ($k_0$) for $\lambda=2$, the first two modes when $\lambda=\phi$, and the first three modes when $\lambda=\sigma$.
    Without forcing, initial conditions are taken such that total energy $E=1$, and first mode, two first modes or three first modes have positive uniform real value, depending on $\lambda=2$, $\phi$ or $\sigma$, while initial smaller scales are zero.
    We observe finite-time blow-up for all three values of $\lambda$ in the two cases.
    Numerical results are plotted in \cref{fig:euler1d}, and scaling exponents are summarized in \cref{tab:exp}.

    \begin{table}
        \caption{\label{tab:exp}Exponents of the inviscid scalings of various quantities measured for the 1D Burgers and the 3D Euler equations within different values of the grid spacing $\lambda$.
        The scalings are with respect to $\tau=1-t/t_b$, where $t_b$ is the blow-up time.
        By definition, the energy spectrum scales like $E(k)\sim k^{-1-2\alpha}$, the maximum value of the vorticity scales like $\omega_{\max}\sim \tau^{-\beta}$, and the width of the singularity strip scales like $\delta\sim\tau^{\mu}$.
        The~$^{(0)}$~superscript indicates a simulation performed with no forcing.
        The~$^*$~superscript indicates a simulation made with a different initial condition.}
        \begin{center}
            \small
            \begin{tabular}{cccccc}
                \br
                & $\lambda$    & $t_b$    & $\alpha$ & $\beta$ & $\mu$  \\
                \mr
                \multirow{5}{*}{1D Burgers} & $2$          & $0.3898$ & $0.37$   & $1$     & $1.55$ \\
                & $\phi$       & $0.5193$ & $0.37$   & $1$     & $1.55$ \\
                & $\sigma$     & $0.4300$ & $0.37$   & $1$     & $1.55$ \\
                & $2^{(0)}$    & $0.2687$ & $0.37$   & $1$     & $1.55$ \\
                & $\phi^{(0)}$ & $0.1460$ & $0.37$   & $1$     & $1.55$ \\
                \mr
                \multirow{4}{*}{3D Euler }  & $2$          & $0.8481$ & $0.67$   & $1$     & $2.81$ \\
                & $\phi$       & $5.8005$ & $0.67$   & $1$     & $2.83$ \\
                & $\phi^*$     & $0.1542$ & $0.67$   & $1$     & $2.82$ \\
                & $\sigma^*$   & $0.8430$ & $0.67$   & $1$     & $2.67$ \\
                \br
            \end{tabular}
        \end{center}
    \end{table}

    The maximum of the gradient $\omega_{\max}(t)=\max_k \vert k u(k)\vert$ blows up in finite time, following the self-similar law
    \begin{equation}
        t_b \omega_{max}\sim\frac{1}{\tau}, \quad \tau = 1-\frac{t}{t_b},
        \label{eq:ssblowup}
    \end{equation}
    displayed in \cref{fig:OmaxAdimBurgers}.
    While the blow-up time depends on the forcing and the value of $\lambda$, the self-similar law~\eqref{eq:ssblowup} is independent of these variables.
    This law also holds for the original continuous model.
    Indeed, differentiating the classical Burgers equation $\partial_t u + u\partial_xu = 0$ with respect to $x$, we get that the space derivative $\omega = -\partial_x u$ obeys $d\omega/dt=\partial_t\omega+u\partial_x\omega=\omega^2$, whose solution is exactly \cref{eq:ssblowup} with $t_b=1/\omega(t=0)$.
    
    We also check that the energy spectrum evolution is universal, in the sense that it only depends on $\tau$.
    This is illustrated in \cref{fig:SpectBurgers}, where spectra for different $\lambda$ but similar $\tau$ are shown to overlap.
    As $\tau$ approaches zero, the energy spectrum gradually widens towards larger values of $k$, developing a power-law $E(k) \sim k^{-2\alpha-1}$ with $1+2\alpha=1.733$, which corresponds to the scaling law $u(k) \sim k^{-\alpha}$ with $\alpha=0.367$.
    Such asymptotics agrees with exact results from the renormalization group formalism applied to the Sabra shell model~\cite{FTBC22}.

    Finally, we compute the analyticity strip width $\delta$ as the solutions approach the blow-up.
    This is done using the formula~\eqref{eq:singustripLL} with $\xi+1=\alpha$.
    The result is shown in \cref{fig:deltaBurgers}.
    We verify that $\delta$ decays to zero in finite time, following a power law $\delta\sim \tau^{\mu}$, with $\mu=1.546$.
    This decay is also universal and does not depend on the value of $\lambda$ or the forcing.
    The width of the analyticity strip is closely associated with $k_{\max}$, defined as the wavenumber at which $\omega$ attains its maximum value.
    Indeed, we see in \cref{fig:deltaAdimBurgers} that $k_{\max}\delta$ is approximately constant in time.
    This is in agreement with the asymptotic \cref{eq:singustripLL}, which implies that $\omega_{\max}$ is achieved at $k_{\max}\sim 1/\delta$.\

    The self-similar law~\eqref{eq:ssblowup} is valid for all $\lambda$ in average only.
    The figures show that the blow-up looks truly self-similar only for the values $\lambda = 2$ and $\lambda = \phi$.
    The oscillations in the case $\lambda = \sigma$ suggest a different blow-up scenario (e.g.\ quasi-periodic or chaotic).
    A detailed analysis of this is left for future work.

    \subsection{Viscous dissipation}
    We now introduce viscous dissipation ($\gamma=1$) and study how the dynamical behaviour depends upon the viscosity parameter $\nu$.
    This section restricts the analysis to the value $\lambda = 2$.
    We introduce a force at the large scale, whose amplitude is adapted dynamically so that the total power input is constant in time ($f_{k=k_0} = P u_{k=k_0} / |u|_{k=k_0}^2$, where $P=1$).

\begin{figure}
        \subfloat[\label{fig:SpectBurgersViscAll}]{%
            \includegraphics[width=0.5\textwidth]{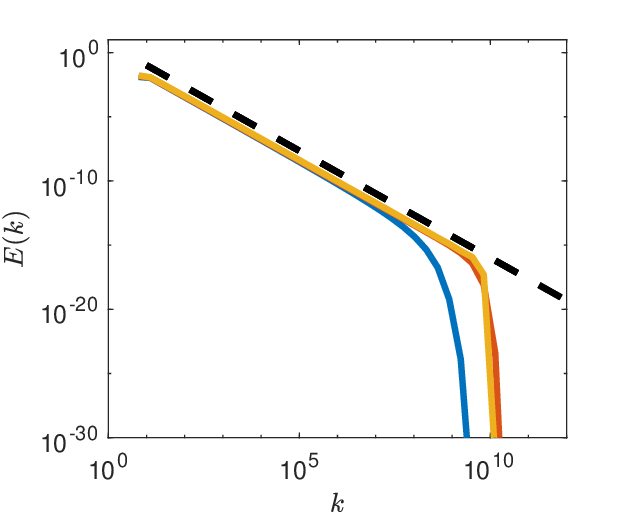}
        }
        \subfloat[\label{fig:OmaxBurgersViscAll}]{%
            \includegraphics[width=0.5\textwidth]{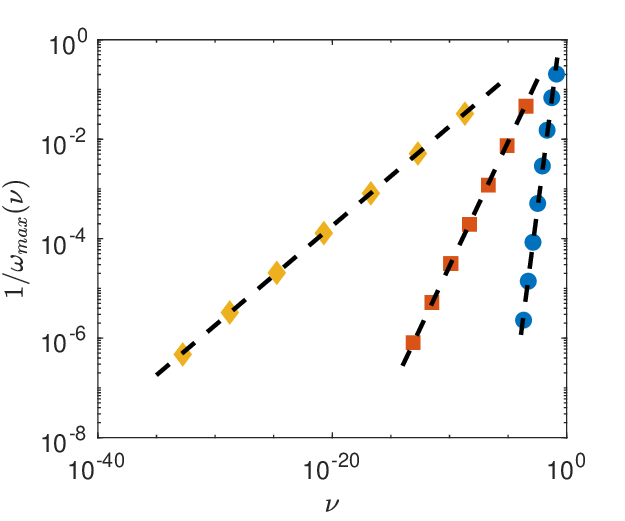}}\\
        \subfloat[\label{fig:deltaBurgersViscAll}]{%
            \includegraphics[width=0.5\textwidth]{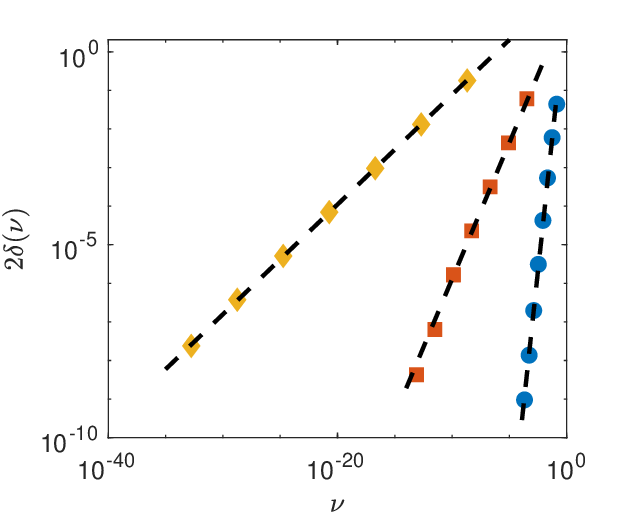}
        }
        \subfloat[\label{fig:deltaAdimBurgersViscAll}]{%
            \includegraphics[width=0.5\textwidth]{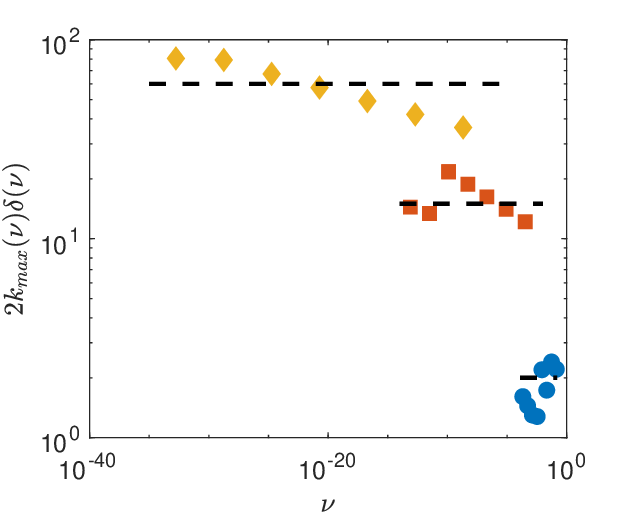}
        }
        ~\caption{Influence of the type of viscosity on the stationary dynamics of the viscous 1D Burgers equation for $\lambda=2$ and $\gamma=0.5$ (hypo-viscous case, blue circle), $\gamma=1$ (viscous case, red squares) and $\gamma=2$ (hyperviscous case, yellow diamond).
            (\ref{fig:SpectBurgersViscAll})~Energy spectrum. The black dotted line has a slope $-5/3$;
            (\ref{fig:OmaxBurgersViscAll})~Maximum value of the derivative $1/ \omega_{\max}$ as a function of viscosity. The black dotted line has a slope given in \cref{tab:exphypviscall} for each case.
            (\ref{fig:deltaBurgersViscAll})~Width of the analyticity strip $2\delta$ as a function of viscosity. The black dotted line has a slope given in \cref{tab:exphypviscall} for each case.
            (\ref{fig:deltaAdimBurgersViscAll})~Renormalized width $k_{\max} \delta$ as a function of viscosity.}
        \label{fig:euler1dviscall}
    \end{figure}

    In this setup, the dissipative term is strong enough to prevent the blow-up, and the solution reaches stationarity.
    The energy spectrum develops a power law in the intermediate scales (called the \textit{inertial range}) followed by an exponential decay at larger $k$ -- see \cref{fig:SpectBurgersViscAll}.
    In the inertial range, $E(k) \propto k^{-5/3}$, corresponding to $u(k)\propto k^{-1/3}$.

    The maximum value of the derivative $\omega_{\max}$ is inversely proportional to the viscosity, following the power law $\omega_{\max}\sim \nu^{-\beta}$ with $\beta=0.5$, as shown in \cref{fig:OmaxBurgersViscAll}.
    This scaling law can be derived when assuming a viscosity-independent \textit{anomalous dissipation} $\epsilon > 0$ in the inviscid limit $\nu \to 0$.
    Under this assumption, we have the balance $\nu\omega^2 \sim \epsilon$, which provides $\omega \sim (\epsilon/\nu)^{1/2}$.

    Accordingly, the width of the analyticity strip does not decline to zero. However, it stabilizes at a finite value that depends on the viscosity -- see \cref{fig:deltaBurgersViscAll} -- and follows the power-law scaling $\delta \sim \nu^\mu$, with exponent $\mu=0.7067$.
    This is smaller than expected from a dimensional argument ``a la Kolmogorov'', in which $\epsilon=\nu u^2/\delta^2$, with $u\sim \delta^{1/3}$, would instead predict $\delta \sim \nu^{3/4}$.
    The strip width follows approximately the scaling $\delta\sim 1/k_{\max}$, as shown in \cref{fig:deltaAdimBurgersViscAll}.

    \subsection{Hyper- and hypo-dissipation}
    We have also studied the influence of the dissipation degree $\gamma$ on the various scaling laws.
    This is summarized in \cref{fig:euler1dviscall,tab:exphypviscall}.
    The slope of the spectrum is insensitive to $\gamma$ and displays a $E(k) \sim k^{-5/3}$ law with no intermittency correction.
    On the other hand, the slopes of both the inverse of the maximum gradient and the singularity width increase in absolute value as $\gamma$ is decreased towards $1/3$.
    We defer the discussion about those results to \cref{subsec:scaling}.
    
     \begin{table}
        \caption{\label{tab:exphypviscall} Scaling exponents of various quantities as a function of $\gamma$ measured for the 1D Burgers and the 3D Navier-Stokes equations with grid spacing $\lambda=2$.
        The scalings are with respect to the viscosity $\nu$. By definition,
            the energy spectrum scales like $E(k)\sim^{-1-2\alpha}$, the maximum value of the vorticity scales like $\omega_{\max}\sim \nu^{-\beta}$ and the width of the analyticity strip
            scales like $\delta\sim\nu^{\mu}$.}
        \begin{center}
            \small
            \begin{tabular}{ccccccc}
                \br
                & \multicolumn{3}{c}{1D Burgers} & \multicolumn{3}{c}{3D Navier-Stokes} \\
                \cmidrule(lr){2-4}\cmidrule(lr){5-7}
                $\gamma$ & $\alpha$ & $\beta$ & $\mu$  & $\alpha$ & $\beta$ & $\mu$  \\
                \mr
                $1/3$    & $-$      & $-$     & $-$    & $2/3$    & $1$     & $2.81$ \\
                $1/2$    & $1/3$    & $1.80$  & $2.78$ & $0.5$    & $0.78$  & $1.53$ \\
                $1$      & $1/3$    & $0.50$  & $0.71$ & $0.40$   & $0.39$  & $0.65$ \\
                $2$      & $1/3$    & $0.20$  & $0.28$ & $0.37$   & $0.19$  & $0.27$ \\
                $8$      & $1/3$    & $0.045$ & $0.06$ & $1/3$    & $0.05$  & $0.06$ \\
                \br
            \end{tabular}
        \end{center}
    \end{table}

    \subsection{Critical dissipation degree \texorpdfstring{$\gamma=1/3$}{}}
    \label{subsec:crit_burgers}
    According to~\cite{Ch08}, there are finite time blow-up solutions for the Burgers equation~\eqref{eq:burgers} with $\lambda = 2$ whenever $\gamma < 1/3$.
    However, the theorems say nothing about the limit case $\gamma = 1/3$.
    For this reason, we call this value as being the \textit{critical dissipation degree}.
    It is natural to ask whether the blow-up might or might not occur in this specific situation.
    Here we consider not only $\lambda = 2$, but also extend this question to the other two lattice parameters.

\begin{figure}
        \subfloat[\label{fig:SpectBurgersCrit}]{%
            \includegraphics[width=0.5\textwidth]{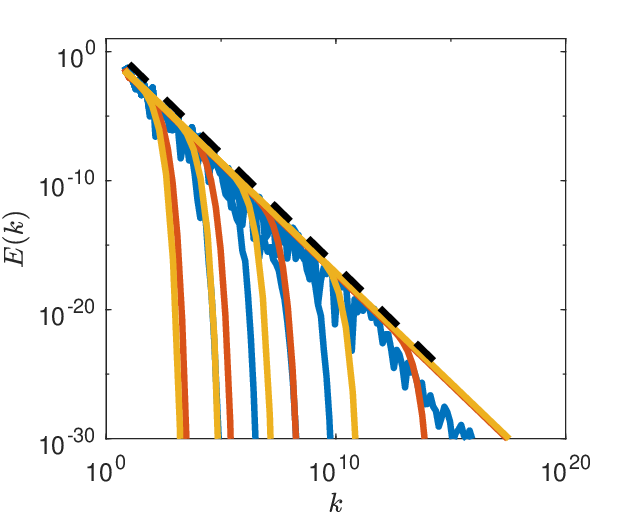}
        }
        \subfloat[\label{fig:OmaxAdimBurgersCrit}]{%
            \includegraphics[width=0.5\textwidth]{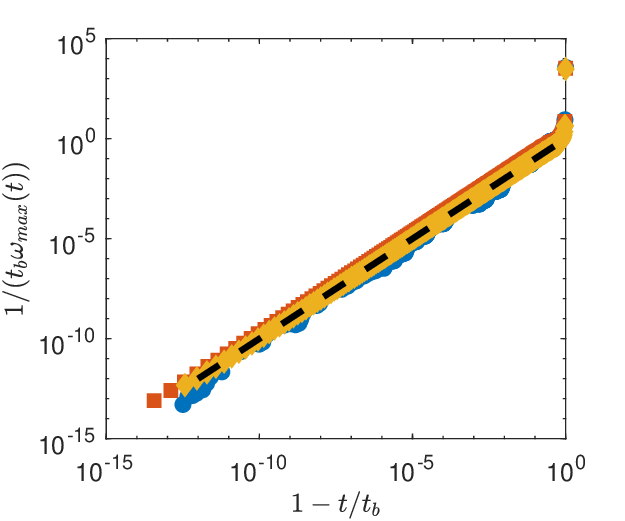}
        }\\
        \subfloat[\label{fig:deltaBurgersCrit}]{%
            \includegraphics[width=0.5\textwidth]{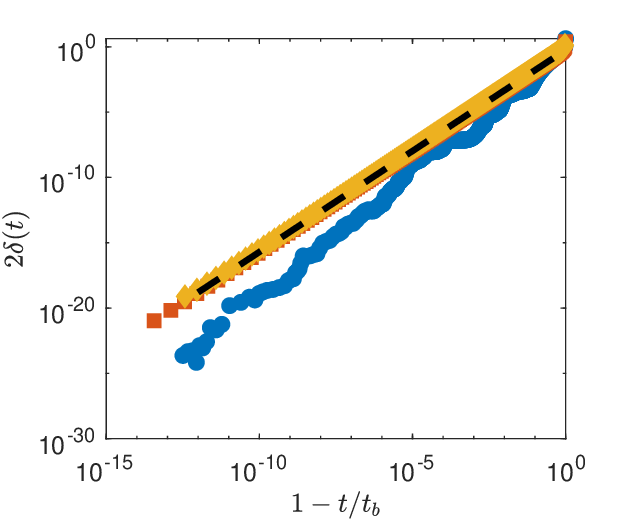}
            \put(-105,30){\includegraphics[width=0.18\textwidth]{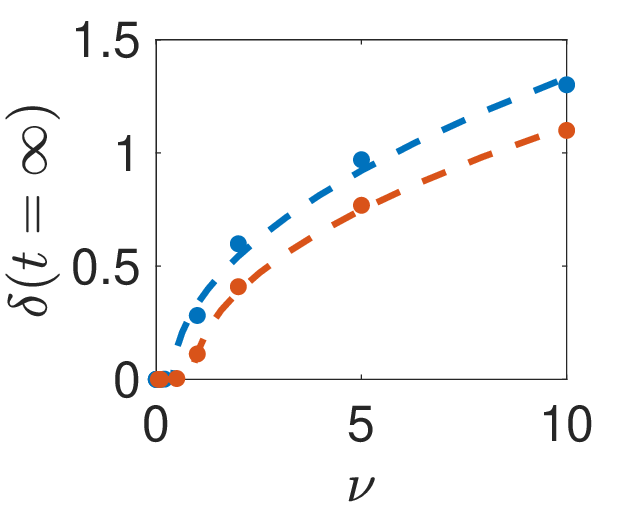}}
        }
        \subfloat[\label{fig:deltaAdimBurgersCrit}]{%
            \includegraphics[width=0.5\textwidth]{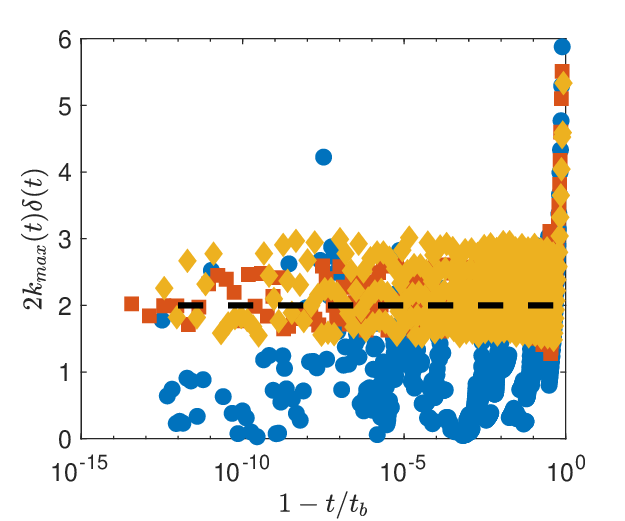}
        }
        ~\caption{Blow-up for the critical ($\gamma=1/3$) 1D Burgers equation, with $\nu=10^{-7}$ and for $\lambda=2$ (yellow), $\lambda=\phi$ (red) and $\lambda=\sigma$ (blue).
            (\ref{fig:SpectBurgersCrit})~Spectra at different renormalized relative time $\tau=1-t/t_b$;
            (\ref{fig:OmaxAdimBurgersCrit})~Maximum value of the derivative $1/t_b \omega_{\max}$ as a function of $\tau$
            (\ref{fig:deltaBurgersCrit})~Width of the analyticity strip $2\delta$ as a function of $\tau$. The insert shows the behaviour of the width of the analyticity strip at $t=\infty$ when the viscosity is increased, for $\gamma=1/3$ (blue data points) (resp. $\gamma=1/4$ (red data points)). The dotted lines are fits of the type $\sqrt{\nu-\nu_c}$, with $\nu_c=0.4$ (resp. $\nu=0.9)$.
            (\ref{fig:deltaAdimBurgersCrit})~Renormalized width $k_{\max} \delta$ as a function of $\tau$. The dotted line has the same scaling and prefactor as in the inviscid blow-up case, see \cref{fig:euler1d}}
        \label{fig:euler1dcrit}
    \end{figure}

    We initialized the flow with the same data as in the inviscid case and set the small viscosity $\nu=10^{-7}$.
    We observed a finite time blow-up for all three $\lambda$, illustrated in \cref{fig:euler1dcrit}.
    The blow-up time is larger than in the inviscid case, but the scaling laws are the same -- both the prefactor and the scaling exponents -- as in the inviscid case.
    The exponents are summarized in \cref{tab:expcrit}.
    The only exception is for the scaling law of $\delta$ in the case $\lambda=\sigma$.
    This might be due to the oscillations in the energy spectrum, making it harder to fit the exponential decreasing, see \cref{fig:SpectBurgersCrit}.

    \begin{table}
        \caption{\label{tab:expcrit}Exponents in the critical case $\gamma=1/3$ of various quantities measured for the 1D Burgers and the 3D Navier-Stokes equations with different values of the grid spacing $\lambda$.
        The scalings are with respect to $\tau=1-t/t_b$, where $t_b$ is the blow-up time. By definition,
            the energy spectrum scales like $E(k)\sim k^{-1-2\alpha}$, the maximum value of the vorticity scales like $\omega_{\max}\sim \tau^{-\beta}$ and the width of the analyticity strip
            scales like $\delta\sim\tau^\mu$.}
        \begin{center}
            \small
            \begin{tabular}{ccccccccc}
                \br
                & \multicolumn{4}{c}{1D Burgers} & \multicolumn{4}{c}{3D Navier-Stokes} \\
                \cmidrule(lr){2-5}\cmidrule(lr){6-9}
                $\lambda$ & $t_b$    & $\alpha$ & $\beta$ & $\mu$  & $t_b$    & $\alpha$ & $\beta$ & $\mu$  \\
                \mr
                $2$       & $0.8497$ & $0.37$   & $1$     & $1.55$ & $7.8194$ & $2/3$    & $1$     & $2.81$ \\
                $\phi$    & $0.5193$ & $0.37$   & $1$     & $1.55$ & $6.51$   & $2/3$    & $1$     & $2.83$ \\
                $\sigma$  & $0.4546$ & $0.37$   & $1$     & $1.84$ & $-$      & $-$      & $-$     & $-$    \\
                \br
            \end{tabular}
        \end{center}
    \end{table}

    In this small-viscosity run, viscosity only delays the blow-up but does not influence the development of the singularity.
    However, we observed a surprising behaviour change when increasing the viscosity to larger values.
    There is a transition between a small-viscosity regime, where finite time blow-up occurs, and a large-viscosity regime, where the blow-up disappears, and the width of the analyticity strip saturates to a finite value -- see insert of \cref{fig:deltaBurgersCrit}.
    The amplitude of $\delta$ seems to follow a critical mean-field behaviour, as it varies like $\delta\sim \sqrt{\nu-\nu_c}$, with $\nu_c\sim 0.4$.
    A similar transition is observed at a lower value of $\gamma$, with $\nu_c$ increasing as $\gamma$ decreases.
    
    This transition is in fact not contradicting the mathematical results by~\cite{Ch08}, since they prove existence of blow-up only for initial conditions larger than a threshold that depends linearly in the viscosity.
    In all our calculations, we start with the same initial conditions.
    This means that for large enough values of viscosity, the initial condition becomes smaller than the threshold, therefore invalidating the hypothesis of the theorem.
    More than that, our numerical results suggest that this hypothesis is actually essential for the result of the theorem and might not be dropped in general.

    \section{3D Euler and Navier-Stokes equations}

    \subsection{Inviscid flow -- Euler equations}

\begin{figure}
        \subfloat[\label{fig:SpectEuler}]{%
            \includegraphics[width=0.5\textwidth]{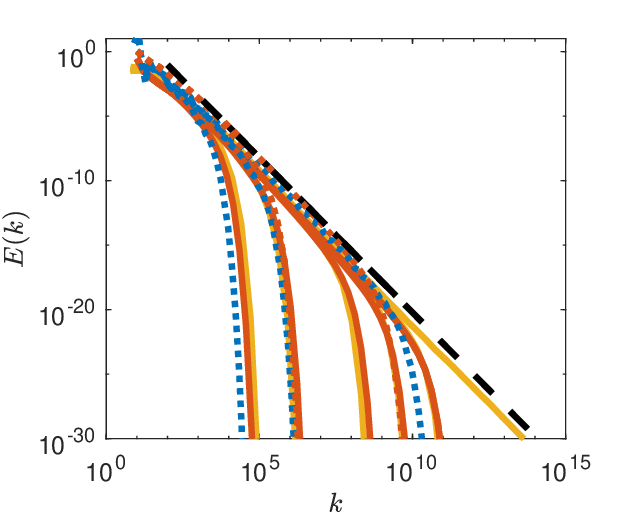}
        }
        \subfloat[\label{fig:OmaxAdimEuler}]{%
            \includegraphics[width=0.5\textwidth]{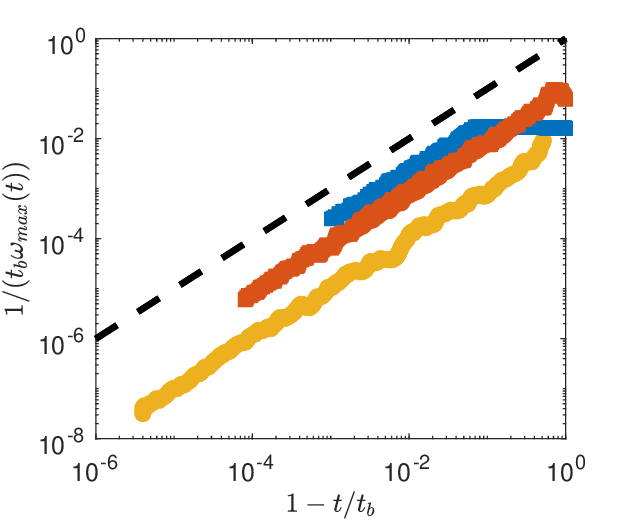}
        }\\
        \subfloat[\label{fig:deltaEuler}]{%
            \includegraphics[width=0.5\textwidth]{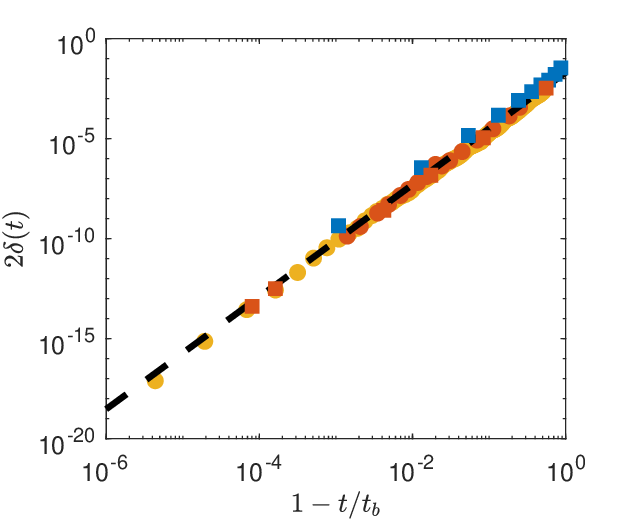}
        }
        \subfloat[\label{fig:deltaAdimEuler}]{%
            \includegraphics[width=0.5\textwidth]{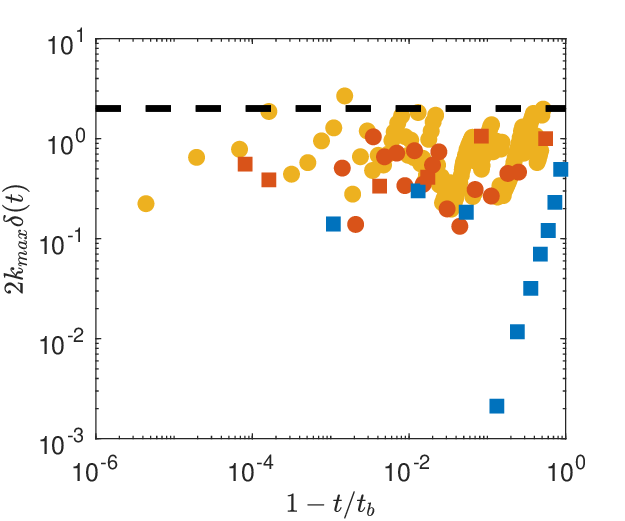}
        }
        ~\caption{Inviscid blow-up for the 3D Euler equations for $\lambda=2$ (yellow), $\lambda=\phi$ (red) and $\lambda=\sigma$ (blue).
            (\ref{fig:SpectEuler})~Spectra at different renormalized relative time $\tau=1-t/t_b$, from $0.2542$ to $0.00001$ from left to right. Spectra with continuous lines and dotted lines correspond to different initial conditions. The black dotted line has a slope of $-7/3$;
            (\ref{fig:OmaxAdimEuler})~Maximum value of derivative $1/t_b \omega_{\max}$ as a function of $\tau$. The black dotted line has a slope of $1$;
            (\ref{fig:deltaEuler})~Width of analyticity strip $2\delta$ as a function of $\tau$; The black dotted line has a slope $2.805$.
            (\ref{fig:deltaAdimEuler})~Renormalized width $k_{\max} \delta$ as a function of $\tau$. In panels~\ref{fig:OmaxAdimEuler}, \ref{fig:deltaEuler} and~\ref{fig:deltaAdimEuler}, we used different symbols for different initial conditions: circles, and squares.}
        \label{fig:euler3d}
    \end{figure}

    We now turn to the full three-dimensional incompressible fluid dynamics on log-lattices, starting with the inviscid Euler equations.
    We consider here the three lattice spacings $\lambda$.
    In order to test universality, we ran the case $\lambda = \phi$ with two different incompressible random  initial conditions, differing by their range of scales. Default initial conditions are defined at large scale $|k| < 3 k_0$, while the other (denoted by a star $*$)  are defined at   scales $|k|<k_0\lambda^3$.
    We observed a finite-time blow-up in all setups, in agreement with previous results documented in~\cite{CM18,CM21}.
    Here, we observe that while the blow-up time depends on the initial conditions, the dynamics become universal when plotted in non-dimensional variables, as illustrated in \cref{fig:euler3d}.
    The spectra for distinct values of $\lambda$ overlap when plotted at the same non-dimensional times $\tau = 1-t/t_b$, as evidenced in \cref{fig:SpectEuler}.
    The slope of the power law in the inertial range is steeper than in 1D Burgers, with a value very close to $-7/3$.
    This is the slope expected for a helicity cascade.
    Our exponent is slightly smaller than those found in some direct numerical simulations of the Euler equations, where a $E(k) \sim k^{-3}$ spectrum is observed~\cite{OPC12,BB12}, but comparable to the value $2.33$ obtained in more recent simulations~\cite{FKMW22}.

    The maximum value of the vorticity $\omega_{\max}$ diverges during the blow-up, as shown in \cref{fig:OmaxAdimEuler}.
    Its asymptotic scaling is the same as for the maximum gradient in the 1D Burgers equation, given by \cref{eq:ssblowup}.
    However, contrarily to the 1D case, the constant in front of the power law varies as a function of $\lambda$ and is not simply given by $1/t_b$.
    This is not too surprising given the 3D nature of the flow, which prevents the application of the simple blow-up argument used for 1D Burgers.
    However, as $\lambda$ is decreased towards~$1$, the non-dimensional curve becomes closer to the exact asymptotic law.

    Approaching the blow-up, the width of the analyticity strip decays to zero with a power law $\delta \sim \tau^\mu$ with exponent $\mu\approx 2.81$ -- see \cref{fig:deltaEuler}.
    This is larger than in 1D Burgers.
    This decay is also universal, as it does not depend on $\lambda$.
    However, it does not show a simple dependence with $k_{\max}$ as seen in \cref{fig:deltaAdimEuler}.
    This might be related to the chaotic nature of the blow-up attractor~\cite{CM21}.

    \subsection{Viscous dissipation -- Navier-Stokes equations}
\begin{figure}
        \subfloat[\label{fig:SpectEulerViscAll}]{%
            \includegraphics[width=0.5\textwidth]{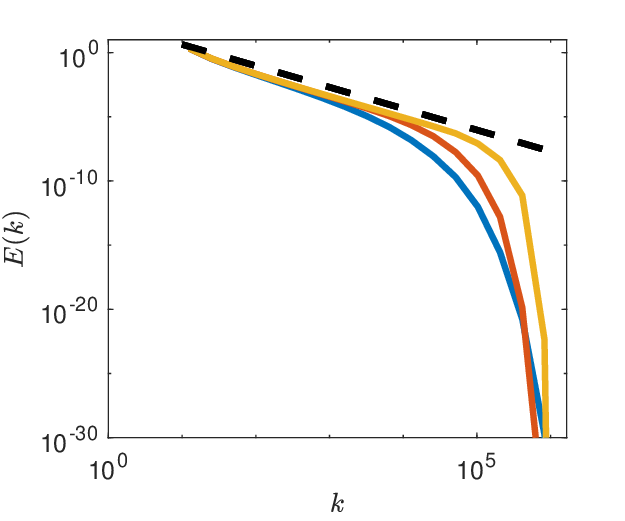}
        }
        \subfloat[\label{fig:OmaxEulerViscAll}]{%
            \includegraphics[width=0.5\textwidth]{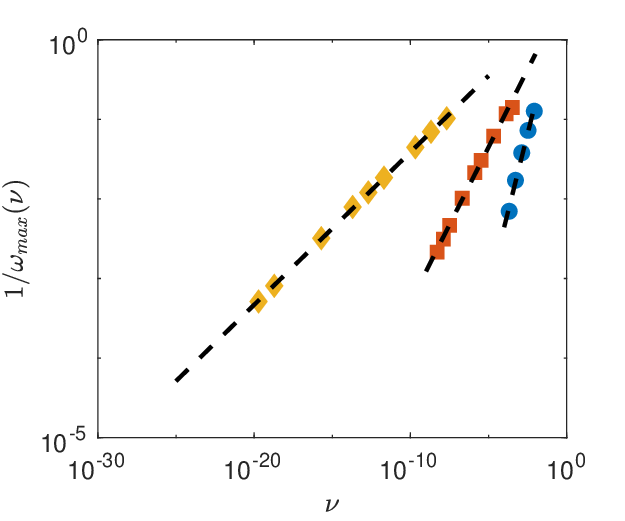}
        }\\
        \subfloat[\label{fig:deltaEulerViscAll}]{%
            \includegraphics[width=0.5\textwidth]{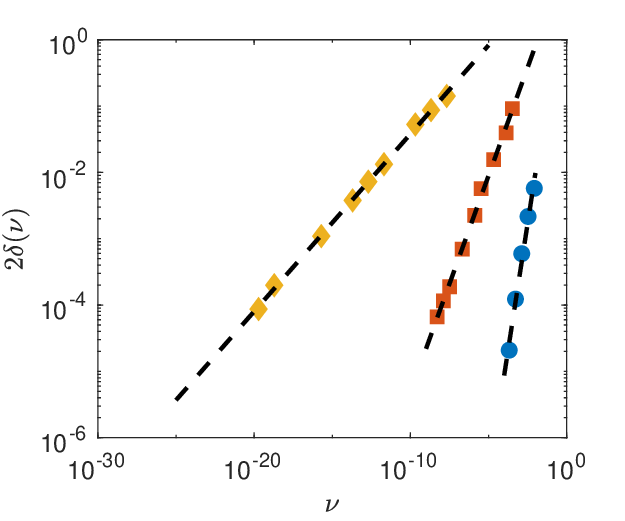}
        }
        \subfloat[\label{fig:deltaAdimEulerViscAll}]{%
            \includegraphics[width=0.5\textwidth]{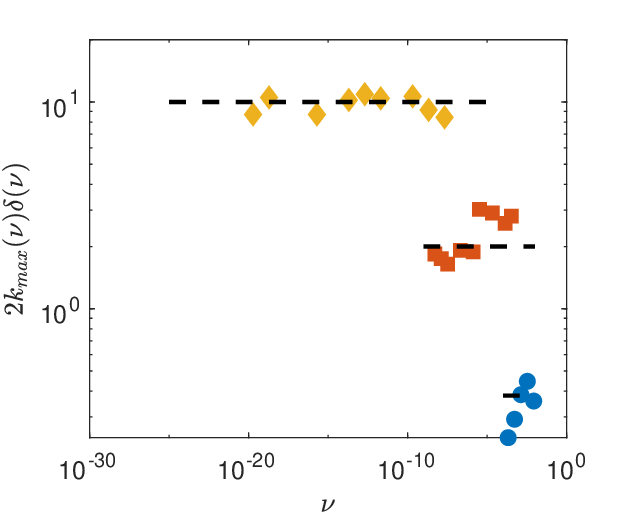}
        }
        ~\caption{Stationary dynamics for the 3D Navier-Stokes equations for $\lambda=2$ and $\gamma=0.5$ (hypo-viscous case, blue circle), $\gamma=1$ (viscous case, red squares) and $\gamma=2$ (hyperviscous case, yellow diamond).
            (\ref{fig:SpectEulerViscAll})~Energy spectrum. The black dotted line has a slope $-5/3$;
            (\ref{fig:OmaxEulerViscAll})~Maximum value of the derivative $1/ \omega_{\max}$ as a function of viscosity. The black dotted line has a slope given in \cref{tab:exphypviscall} for each case.
            (\ref{fig:deltaEulerViscAll})~Width of the analyticity strip $2\delta$ as a function of viscosity. The black dotted line has a slope given in \cref{tab:exphypviscall} for each case.
            (\ref{fig:deltaAdimEulerViscAll})~Renormalized width
            $k_{\max} \delta$ as a function of viscosity. }
        \label{fig:euler3dviscall}
    \end{figure}

    We now add the viscous term with $\gamma = 1$ and a constant-power forcing.
    The solutions achieve a statistically stationary state in this framework, whose average scalings are depicted in \cref{fig:euler3dviscall}.
    Like in 1D Burgers, the energy spectra display a power law until the solution reaches the viscous scale, with an inertial range widening as $\nu$ decreases.
    The slope of the energy spectrum is slightly steeper than Kolmogorov's $-5/3$, with an intermittency correction of around $0.13$.
    Accordingly, the scaling exponent $\alpha$ for $u(k) \sim k^{-\alpha}$ is $\alpha = 0.40$, slightly larger than $1/3$.
    The maximum vorticity $\omega_{\max}$ increases with decreasing viscosity, following the power law $\omega_{\max} \sim \nu^{-\beta}$ with an exponent $\beta = 0.39$ lower than Kolmogorov's $1/2$.

    The width of the analyticity strip decays with viscosity as $\delta \sim \nu^\mu$ with an exponent $\mu = 0.65$ -- see \cref{fig:deltaEulerViscAll}.
    Such decay is less intense than in the 1D Burgers equation.
    Nevertheless, the dependence of $\delta$ on $1/k_{\max}$ in the Navier-Stokes case is sharper, as one verifies by comparing \cref{fig:deltaAdimEulerViscAll} with \cref{fig:deltaAdimBurgersViscAll}.


    \subsection{Hyperviscous dissipation}
    We now consider what happens in the hyperviscous case $\gamma > 1$.
    We keep the constant-power forcing to reach stationary states.

    For $\gamma = 2$, we still observe a power-law
    energy spectrum followed by an exponential cut-off at the viscous scales -- see \cref{fig:SpectEulerViscAll}.
    The inertial range keeps widening as $\nu$ is decreased.
    The slope of the energy spectrum is very close, but slightly steeper than  $-5/3$. The exact fitting provides us an intermittency correction around $0.07$, corresponding to $\alpha=0.37$, see table~\ref{tab:exphypviscall}.
    The maximum vorticity $\omega_{\max}$ increases with decreasing viscosity like a power law, with an exponent $\beta=0.19$ lower than usual ($\gamma = 1$) viscous case.
    The width of the analyticity strip decays with viscosity with an exponent $\mu=0.26$.
    Like in the viscous case, $\delta$ appears to scale simply like $1/k_{\max}$, as seen on \cref{fig:deltaAdimEulerViscAll}.

    The above results suggest that the intermittency corrections in the energy spectra are smaller for hyperdissipation.
    Indeed, as the dissipation degree $\gamma$ increases, the exponent $\alpha$ converges towards Kolmogorov's $1/3$, see table~\ref{tab:exphypviscall}.
    We checked that for the stronger degree $\gamma=8$, they vanish completely, and the dependence of $\omega_{\max}$ and $\delta$ on $\nu$ become very weak.
    This is explained by the very sharp viscous cut-off due to the hyperviscous dissipation.
    Indeed, the equivalent of the Kolmogorov scale $k_d$ in the hyperviscous case relates to $\nu$ as $k_d \sim \nu^{1/(1-1/3-2\gamma)}$, becoming independent of viscosity in the limit $\gamma \to \infty$.
    For $\gamma=2$ the dependence is $\delta\sim k_d^{-1} \sim \nu^{0.3}$, close to what is observed for the scaling of the singularity strip width.

    \subsection{Hypoviscous dissipation}
    The case with hypoviscous dissipation $1/3<\gamma<1$ is qualitatively similar to the viscous and hypervisous cases -- see \cref{fig:euler3dviscall}.
    Exponents, however, are steeper.
    The corresponding values are reported in \cref{tab:exphypviscall}.
    The energy spectrum develops a slope corresponding to the exponent $\alpha = 0.5$, which is steeper than Kolmogorov's $1/3$ but milder than Euler's $2/3$ on log-lattices.
    The singularity width appears again to be controlled by the wave number corresponding to the maximum vorticity -- see \cref{fig:deltaAdimEulerViscAll}.
    On the other hand, the maximum vorticity grows much more rapidly than in the viscous case, with an exponent twice as big, as shown in \cref{fig:OmaxEulerViscAll}.
    This may indicate that we are approaching a critical dissipation degree, below which finite-time blow-up will occur.

    \subsection{Critical dissipation degree \texorpdfstring{$\gamma=1/3$}{}}
    \label{subsec:ns_crit}
    The asymptotics of Kolmogorov's length scale for a flow with a general dissipation degree predicts the breakdown of the viscous cut-off when $\gamma$ approaches the critical value $1/3$.
    Indeed, the dissipation scale $k_d$ is obtained from the dimensional balance between the convective and the dissipative terms $k_du_d^2 \sim \nu k_d^{2\gamma}u_d$.
    On the other hand, Kolmogorov's theory states that $u_d \sim \epsilon^{1/3}k_d^{-1/3}$ for the energy dissipation rate $\epsilon$, which has a finite positive value in the inviscid limit.
    Together, these expressions yield
    \begin{equation}
        \label{eq:k_d_scaling}
        k_d \sim \epsilon^{\frac{1}{6(\gamma - 1/3)}}\nu^{\frac{1}{2(1/3-\gamma)}},
    \end{equation}
    which, for sufficiently small $\nu$, provides $k_d \to +\infty$ when $\gamma \searrow 1/3$.
    For this reason, we call $\gamma = 1/3$ the \textit{critical dissipation exponent}, the value at which we expect that the dissipative term is no longer strong enough to prevent a finite-time singularity.
    We recall this was the case for the 1D Burgers equation on log-lattices.

    Motivated by the above arguments, we investigate the critical hypo-diffusive degree in the full 3D system on log-lattices.
    The following analysis considers the spacings $\lambda = 2$ and $\lambda = \phi$.
    The initial data is the same as we used in the inviscid simulations, and viscosity is the same $\nu = 10^{-7}$. 

\begin{figure}
        \subfloat[\label{fig:SpectEulerCrit}]{%
            \includegraphics[width=0.5\textwidth]{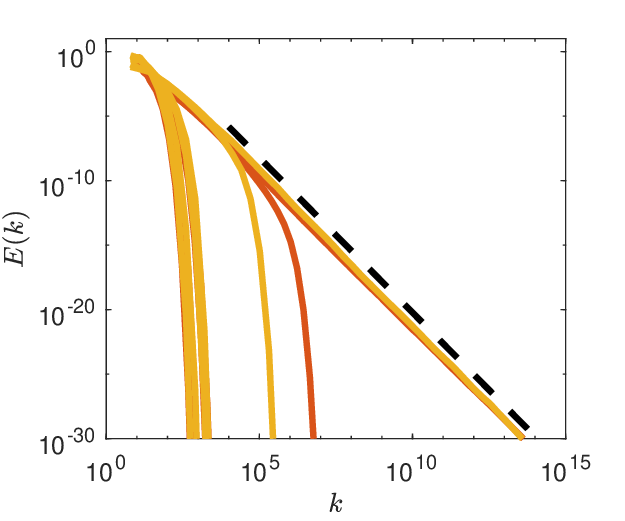}
        }
        \subfloat[\label{fig:OmaxAdimEulerCrit}]{%
            \includegraphics[width=0.5\textwidth]{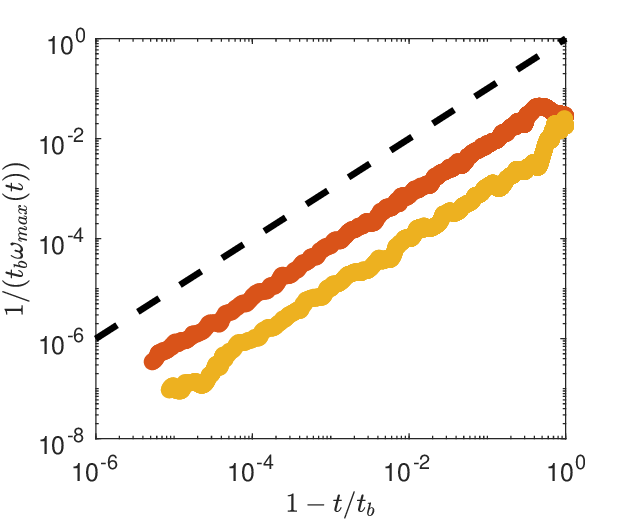}
        }\\
        \subfloat[\label{fig:deltaEulerCrit}]{%
            \includegraphics[width=0.5\textwidth]{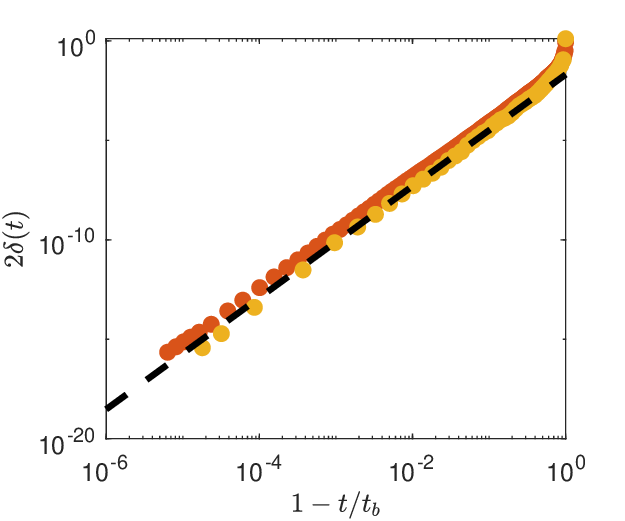}
        }
        \subfloat[\label{fig:deltaAdimEulerCrit}]{%
            \includegraphics[width=0.5\textwidth]{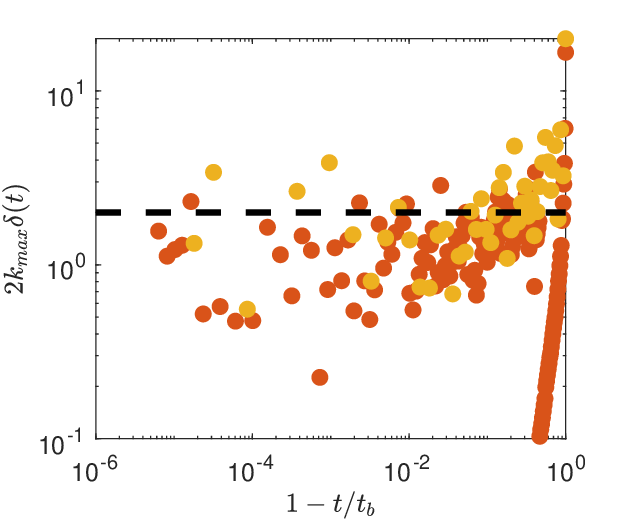}
        }
        ~\caption{Blow-up for the critical ($\gamma = 1/3$) 3D Navier-Stokes equations, for $\lambda=2$ (yellow) and $\lambda=\phi$ (red).
            (\ref{fig:SpectEulerCrit})~Spectra at different renormalized relative time $\tau=1-t/t_b$;
            (\ref{fig:OmaxAdimEulerCrit})~Maximum value of the derivative $1/t_b \omega_{\max}$ as a function of $\tau$
            (\ref{fig:deltaEulerCrit})~Width of the analyticity strip $2\delta$ as a function of $\tau$
            (\ref{fig:deltaAdimEulerCrit})~Renormalized width $k_{\max} \delta$ as a function of $\tau$. The dotted line has the same scaling and prefactor as in the inviscid blow-up case, see \cref{fig:euler3d}}
        \label{fig:euler3dcrit}
    \end{figure}

    In this regime, we observed a finite time blow-up for the two values of $\lambda$, illustrated in \cref{fig:euler3dcrit}.
    Like in 1D Burgers, the blow-up time is larger than in the inviscid case, but the scaling laws are the same.
    This is summarized in \cref{tab:expcrit}.
    The slope of the energy spectrum remains $-7/3$.
    For $\lambda = 2$ and $\nu = 10^{-3}$, the dynamics becomes stationary, meaning there is as in 1D a phase transition, but between $\nu = 10^{-3}$ and $10^{-7}$, smaller than $\nu_c \sim 0.4$ in 1D.

    \section{Discussion}

    \subsection{Scaling laws}
    \label{subsec:scaling}
\begin{figure}
        \subfloat[\label{fig:alpdegamma}]{%
            \includegraphics[width=0.33\textwidth]{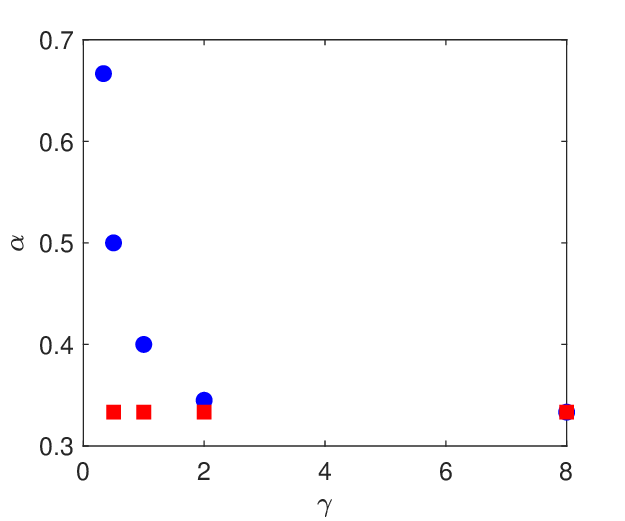}
        }
        \subfloat[\label{fig:betadegamma}]{%
            \includegraphics[width=0.33\textwidth]{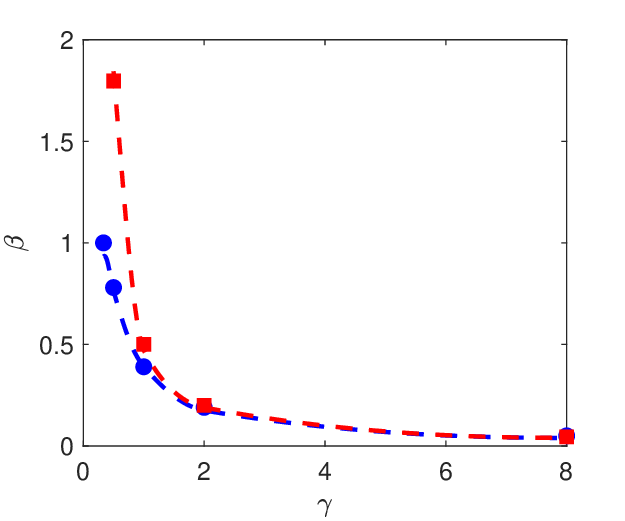}
        }
        \subfloat[\label{fig:mudegamma}]{%
            \includegraphics[width=0.33\textwidth]{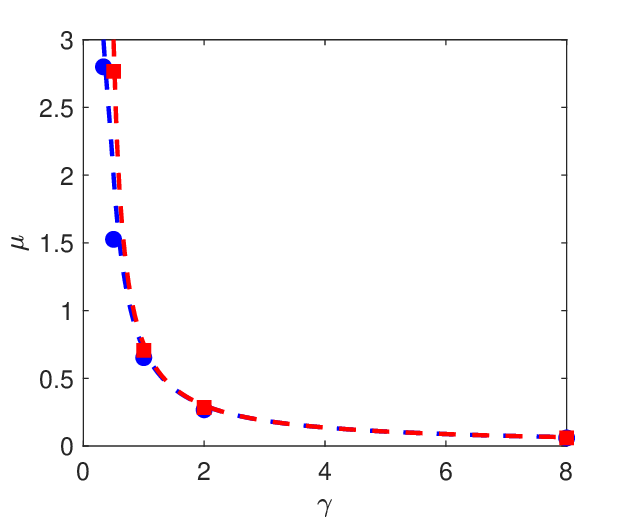}
        }
        ~\caption{Variation of the scaling exponents as a function of the dissipation degree $\gamma$ for the 1D Burgers (blue circle) and the 3D Navier-Stokes (red square) equations with $\lambda=2$.
            (\ref{fig:alpdegamma})~For the scaling of the velocity $u\sim k^{-\alpha}$.
            (\ref{fig:betadegamma})~For the scaling of the maximum vorticity $\omega_{\max}\sim \nu^{ -\beta}$. The data points are reported from \cref{tab:expcrit,tab:exphypviscall}, while the dotted lines correspond to \cref{eq:linkexponents};
            (\ref{fig:mudegamma})~For the scaling of the singularity strip $\delta\sim \nu^{\mu}$. The data points are reported from \cref{tab:expcrit,tab:exphypviscall}. while the dotted lines correspond to \cref{eq:equamu}.}
        \label{fig:scalinglaws}
    \end{figure}

    The variations of the scaling exponents with respect to the diffusion exponent $\gamma$ are shown in \cref{fig:scalinglaws}.

    Predictions for the scaling laws are possible using simple dimensional arguments if we impose $\delta\sim 1/k_{\max}$, as empirically observed.
    Indeed, from $u\sim k^{-\alpha}$ and $\omega\sim k u$, we get $\omega_{\max}\sim k_{\max}^{1-\alpha}\sim \delta^{\alpha-1}$ so that we get:
    \begin{equation}
        \beta=\mu(1-\alpha).
        \label{eq:linkexponents}
    \end{equation}
    This fixes a link between the 3 exponents that is well satisfied -- see \cref{fig:betadegamma}.
    On the other hand, one can connect $\mu$ and $\alpha$ by extending the argument fixing the Kolmogorov scale to hypo and hyper-viscous cases: we impose that $k_{\max}$ is fixed by the condition that the viscous term balances the non-linear term $\nu k_{\max}^{2\gamma} u_{\max}\sim k_{\max} u_{\max}^2$.
    Using $u_{\max}\sim k_{\max}^{-\alpha}$ and $\delta\sim 1/k_{\max}$ we then get:
    \begin{equation}
        \mu=-\frac{1}{1-\alpha-2\gamma},
        \label{eq:equamu}
    \end{equation}
    This prediction is tested in \cref{fig:mudegamma} and is well satisfied.
    Without loss of generality, the only free parameter can be taken as  $\alpha(\gamma)$.
    In the limit $\gamma\to 1/3$, we can fix it by imposing that
    $\beta=1$, which is the scaling corresponding to conservation of the circulation of $u$~\cite{Po18}.
    From \cref{eq:equamu,eq:linkexponents}, we then get
    $\alpha=1-\gamma=2/3$, corresponding to a helicity cascade.
    In all other cases, we have no clear theories to predict the variations of $\alpha$ with $\gamma$.
    Notably, when $\gamma\to \infty$, we recover $\alpha=1/3$ corresponding to an energy cascade.

    \subsection{Interest of the critical case}
    The critical case $\gamma=1/3$ is more than purely academic: renormalization group (RNG) analysis of NSE in Fourier space~\cite{YO86} indeed shows that the fixed point of the equations corresponds to a Navier-Stokes equation with turbulent viscosity scaling like $A\epsilon^{1/3}k^{-4/3}$, where $A$ is a constant with value $A=0.1447$ in 1D and $A=0.4926$ in 3D. This corresponds exactly to \cref{eq:loglatice}, with $\gamma=1/3$ and $\nu=A\epsilon^{1/3}$.
    This model is sometimes used as a subgrid model of turbulence~\cite{LDN01}.
    In that respect, it is interesting that the transition viscosity found in \cref{subsec:crit_burgers,subsec:ns_crit} (at constant injected power, i.e.\ $\epsilon=1$) is very close to the RNG value in 1D\@.
    On the one hand, this guarantees that the size of the inertial range is very wide, in agreement with the RNG picture of scale invariant solutions.
    On the other hand, this means that the solution is very close to a blow-up, which could have implications regarding the stability of this subgrid scheme.

    \subsection{Implications for real Euler or Navier-Stokes?}
    The log-lattices simulations we performed cannot be seen as an exact model of the Euler or Navier-Stokes equations because they remove by construction many non-linear interactions of the original equations, especially the non-local one.
    Nevertheless, because they obey the same conservation laws and symmetries, they may capture some scaling laws of the original equation more accurately.
    Comparing our findings with the few results on the topic is engaging.

    Regarding the Euler equation, recent high-resolution numerical simulation in the axisymmetric case by~\cite{KSP22} explored the scaling of the singularity strip in the blowing situation proposed by~\cite{LH14}.
    They found an exponent $\mu=2.6\pm 0.5$, which is compatible with the value $2.8\pm0.1$ that we get from \cref{tab:exphypviscall}.
    Unfortunately, they do not provide an estimate of the slope of the energy spectrum.
    Previous older results in the Taylor-Green vortex~\cite{BMVP+92,BB12} found a steeper spectrum corresponding $\alpha\sim 1$.
    However, spectra with exponent matching our $-7/3$ value were observed in the early stage of recent simulations at larger resolution~\cite{FKMW22}.
    Therefore, the main characteristics of blow-up in log-lattices simulations agree with the most recent results observed in the traditional DNS of the Euler equation.

    Regarding the Navier-Stokes equations, we can look at two recent results.
    The first one by~\cite{Du22} finds a value of $\mu=0.89$ using recent DNS of NSE\@.
    This value is larger than the value we found in the present paper, corresponding to $\mu=0.65$.
    Another recent result~\cite{BPBY19} estimates $\beta$ in 3D NSE. They indeed found that the tail of the PDF of enstrophy scales like $\nu^{0.77}\tau_K^{-2}$, where $\tau_K\sim\nu^{1/2}$ is the Kolmogorov time.
    Identifying such extreme events of enstrophy with $\omega_{\max}^2$, we thus get $\beta_{DNS}\sim 0.88$, which is also much larger than the value we observe in log-lattices $\beta_{LL}\sim 0.39$.
    Note, however, that both DNS values are compatible with \cref{eq:linkexponents,eq:equamu}, provided we choose $\alpha\sim 0$, hinting at the presence of multifractality.
    Log-lattices simulations are generally much less intermittent than DNS~\cite{CM21}, with one dominating exponent (monofractal behaviour).
    Some time ago,~\cite{LDN01} linked the intermittency properties of NSE with non-local interactions, which is coherent with this observation.
    Therefore, the difference between log-lattices simulations and DNS could be explained by differences in the amount of non-local interactions.

    \section*{Acknowledgments}
    This work received funding from the Ecole Polytechnique, from ANR TILT grant agreement no.
    ANR-20-CE30-0035, and from ANR BANG grant agreement no.
    ANR-22-CE30-0025.
    CC is thankful for the financial support from CEA during his visit to Paris-Saclay, where this work was elaborated.
    AM was supported by CNPq grant 308721/2021-7 and FAPERJ grants E-26/201.054/2022, E-26/210.874/2014.

    \section*{References}
    \bibliographystyle{unsrt.bst}
    \bibliography{fractalsLL}

\end{document}